\documentclass{article}
\usepackage{amsmath,amssymb,amsfonts}
\usepackage{graphicx}
\usepackage{color}

\def\deg {\mathrm{deg}}

\def\log {\mathrm{log}}
\def\GCD {\mathrm{GCD}}
\def\sn{\mathrm{sn}}

\newtheorem{theorem}{Theorem}
\newtheorem{proposition}{Proposition}
\newtheorem{corollary}{Corollary}
\newtheorem{lemma}{Lemma}

\newtheorem{example}{Example}

\newenvironment{proof}[1][Proof]{\noindent\textit{#1.} }{\hfill$\Box$\medskip}

\title{Marden theorem and Poncelet-Darboux curves }
\author{Vladimir Dragovi\' c}

\date{}

\begin{document}

\maketitle

\medskip

\centerline{Mathematical Institute SANU}

\centerline{Kneza Mihaila 36, 11000 Belgrade, Serbia}

\smallskip

\centerline{Mathematical Physics Group, University of Lisbon}

\smallskip

\centerline{e-mail: {\tt vladad@mi.sanu.ac.rs}}

\

\begin{abstract}

\smallskip

The Marden theorem of geometry of polynomials and the great
Poncelet theorem from projective geometry of conics by their
classical beauty occupy very special places. Our main aim is to
present a strong and unexpected relationship between the two
theorems. We establish a dynamical equivalence between the full
Marden theorem and the Poncelet-Darboux theorem. By introducing a
class of {\it isofocal deformations}, we construct  morphisms
between the Marden curves and the Poncelet-Darboux curves. Then we
present effective criterion in terms of pair of polynomials which
defines a Poncelet-Darboux curve of degree $n-1$, for complete
decomposition of the curve  on $(n-1)/2$ conics if $n$ is odd; if
$n$ is even, complete decomposition consists of $(n-2)/2$ conics
and a line. This is an important question in the study of special,
'tHooft, instanton bundles.
\end{abstract}

\newpage

\tableofcontents

\newpage

\section{Introduction}\label{sec:intro}

\medskip
Although it is impossible to distinguish the nicest mathematical
result, probably many mathematicians would agree that the Marden
theorem of geometry  of polynomials and the great Poncelet theorem
from projective geometry of conics by their classical beauty
occupy very special places. Our main aim is to present a strong
and unexpected relationship between the two theorems. We establish
a dynamical equivalence between the full Marden theorem and the
Poncelet-Darboux theorem. By introducing a class of {\it isofocal
deformations}, we construct a morphism between a Marden curve and
a Poncelet-Darboux curve, using the Moser trick and the Flashka
coordinates. As a byproduct, we get complete description of
cyclic-symmetric $n$-correspondences of $\mathbb P^1$. Then we
present effective criterion for complete decomposition of a
transversal Poncelet-Darboux curve of degree $n-1$ on $(n-1)/2$
conics if $n$ is odd and on $(n-2)/2$ conics and a line if $n$ is
even in terms of pair of polynomials which defines the
Poncelet-Darboux curve. This is an important question in the study
of special or 'tHooft instanton bundles (see \cite{Tr},\cite{Tr1},
\cite {HN}, \cite {Ha}, \cite {NT} and references therein). We
finish by introducing a new class of {\it bifocal deformations}
and relating it again to classical Darboux results.

The idea to relate Marden's theorem with Poncelet's theorem
appeared after our careful analysis of the proof of the Siebeck
theorem for $n=3$. For the reader's sake and for self-completion
we included an Appendix with explanation of basic steps of the
proof of the Siebeck theorem, following \cite {Ka} and references
therein.

In the next Section \ref{sec:prel} we formulate theorem of
Siebeck, its generalization - the Marden theorem and theorems of
Poncelet and Darboux.

The Section \ref{sec:isofocal} is devoted to introduction of
isofocal dynamics, which appears to be integrable. By application
of the Moser trick and the Flashka coordinates, we get explicit
trivialization of the dynamics, see Theorem \ref{th:flashka}. A
birational morphism between the data which define a Marden curve
and the data defining a Poncelet-Darboux curve is established in
the Theorem \ref{th:mardendarboux}.

The main question treated in the Section \ref{sec:decomposition}
is about necessary and sufficient conditions on pair of
polynomials which defined a transversal Poncelet-Darboux curve $S$
of degree $n-1$ to be completely decomposed on $(n-1)/2$ conics
for $n$ odd; for $n$ even the complete decomposition assumes
decomposition on $(n-2)/2$ conics and a line. In a generally very
nice paper \cite{Tr} this question has been treated explicitly,
but in unsatisfactory way (See Sec. 5.5 of \cite{Tr}:unfortunately
the conditions mentioned there are neither sufficient nor
necessary). Here, we propose a systematic approach to this
important and difficult question. Very strong necessary conditions
are formulated in the Theorem \ref{th:condition}. It leads to the
study of elliptic coverings of elliptic curves and to the theory
of transformations of elliptic functions which was established by
Jacobi (see \cite{Jac2}). In order to get sufficient conditions,
one needs to examine if above transformation is cyclic, see
Theorem \ref{th:convereseconditions}. For $n=2^km$ with $m$ odd
all above transformations are cyclic for $k=0, 1$, but not for
$k>1$, see Theorem \ref{th:cyclic}.

In Section \ref{subsec:choreography} we demonstrate effectiveness
of the previous considerations. For $n=3, 5, 7$ we give complete
list of pairs of polynomials which define completely decomposable
Poncelet-Darboux curves of degree $n-1$. Moreover, we describe
corresponding initial conditions of the isofocal deformations.

In the last Section \ref{sec:conclusion} we briefly mention three
related problems: the bifocal transformations; the Toma -
Trautmann case of conic component and the positivity problem in
isofocal dynamics, connected with infrapolynomial interpretation
of the Marden theorem.

\medskip

\section{Preliminaries}\label{sec:prel}

\medskip

\subsection{The Marden theorem}\label{subsec:marden}

\medskip

One of the basic theorems in geometric theory of polynomials and
rational functions has a long history and is usually referred as
{\it Marden's} after appearance of the book \cite {Ma}. The
earliest version of this theorem, up to our best knowledge, goes
back to 1864 when Siebeck  (see \cite {Si}) formulated and proved
it for the case of polynomials with simple roots:

\medskip

\begin{theorem} [Siebeck \cite {Si}] \label{th:Siebeck}
Let $P(z)$ be a polynomial of
degree $n\ge 3$ with complex coefficients, such that the zeros
$\alpha_1, \dots,\alpha_n$ are simple and every three
noncollinear. There exists a curve $C$ of class $n-1$ tangent to
every line segment $[\alpha_i, \alpha_j]$ at the midpoint. The
foci of the curve $C$ are zeros of the derivative polynomial
$DP(z)=P'(z)$.
\end{theorem}

\medskip

In the simplest case $n=3$  the curve $C$ is a conic, inscribed in
the triangle formed by the zeros of a polynomial of degree 3 and
tangent to the sides of the triangle at their midpoints. Even in
this, simplest case, the result of the Siebeck theorem is
nontrivial and interesting and attracted lot of attention not only
in the past but also nowadays (see for example \cite {Ka}).

\medskip

Previous results were extended to the cases with not all roots
being simple. Consider a function of the form
\begin{equation}\label{eq:polynomial}
P(z)=(z-\alpha_1)^{m_1}
(z-\alpha_2)^{m_2}\cdots(z-\alpha_n)^{m_n},
\end{equation}
where all $\alpha_k$ are distinct, every three noncollinear and
$N=m_1+m_2+\dots m_n$. The zeros of the derivative $DP$ are
divided into two groups. In the first group are those $\alpha_i$
for which  $m_i>1$. The second group is formed from simple zeros
of the function, in other words from the zeros  of the logarithmic
derivative $LP$ of $P$. Since the positions of the zeros of the
first group are known from the beginning, the interesting part is
location of the members of the second group.

Thus, consider the function $F(z)=LP(z)=d[\log P(z)]/dz$:
\begin{equation}\label{eq:function}
F(z)=\frac{m_1}{z-\alpha_1}+\dots +\frac{m_n}{z-\alpha_n}.
\end{equation}

Historically, the constants $m_i$ in the last expressions were
firstly considered as positive integers. Then, step by step, that
condition has been relaxed up to the condition that $m_i$ are
nonzero real numbers, as we can find in Marden's formulation (see
\cite {Ma}, p. 11, Th. 4.2):

\medskip

\begin{theorem}[\cite {Ma}, p. 11]\label{th:marden}
The zeros of the function
\begin{equation}
F(z)=\frac{m_1}{z-\alpha_1}+\dots +\frac{m_n}{z-\alpha_n}.
\end{equation}
where $m_i$ are real nonzero constants, are the foci of the curve
of class $n-1$ which touches each line-segment $[\alpha_i,
\alpha_j]$ in a point dividing the line segment in the ratio
$m_i:m_j$.
\end{theorem}

\medskip
A good account of a century-long path from the Theorem
\ref{th:Siebeck} to the Theorem \ref{th:marden} one can find in
Marden's book \cite {Ma}, together with references. Thus we are
going to omit them here.

\medskip

\subsection{Darboux Theorem and Poncelet-Darboux curves}

\medskip

Now we pass to quite different subject, which appeared in the
theory of conics in the context of the great Poncelet theorem
\cite {Pon}. An overview of the history of the subject together
with presentation of contemporary state of art one may find in
\cite {DR}, \cite {DR2}. We start here with one of Darboux's
original formulations of his theorem.

\medskip

\begin{theorem}[Darboux, \cite{Dar1} p. 248]\label{th:darboux1}
Si une courbe d'ordre $n-1$ contient tous les points
d'intersection de $n$ tangents a une conique, elle contient aussi
les points d'intersection d'une infinit\'e d'autres systemes de
$n$ tangentes a la m\^eme conique. Chacun de ces systemes est
d\'efini par l'equation
\begin{equation}\label{eq:pencildarboux}
\phi(\rho)+kf(\rho)=0
\end{equation}
ou $k$ d\'esigne une constante arbitraire.
\end{theorem}

\medskip

Darboux had been interested in this matter for about fifty years
and he published several variations of the last theorem (see for
example \cite {Dar3}). Slightly changing terminology from \cite
{Tr1}, we will say that a curve $S$ of degree $n-1$  {\it is
Poncelet-Darboux related to a conic $K$} if the curve $S$ and the
conic $K$ satisfy conditions of the previous theorem. The set of
all such curves of degree $n-1$ which are Poncelet-Darboux related
to a fixed conic $K$ will be denoted as
$$
Pon-Dar_{n-1}(K).
$$

We will say that  the $n$ tangents $t_1, t_2, \dots, t_n$ of a
conic $K$ from the previous theorem {\it form a Poncelet
$n$-polygon} $P_n=T_1T_2\dots T_n$, where $T_i=t_i\cap t_{i+1}$
for $i=1,\dots, n-1$ and $T_n=t_{n}\cap t_1$ if there exists a
conic $C$ such that $T_i\in C$ for $i=1,\dots, n$. In that case we
will say that conics $C$ and $K$ are {\it $n$-Poncelet related.}

Here is a formulation of  the Poncelet Theorem:

\medskip

\begin{theorem}[Poncelet, \cite {Pon}]\label{th:poncelet}
If two conics $C$ and $K$ are $n$-Poncelet related, then there are
infinitely many $n$-polygons circomscribed about $S$ and inscribed
in $C$. Moreover, arbitrary point of the conic $C$ may be chosen
for a vertex of a such Poncelet $n$-polygon.
\end{theorem}

\medskip

If the $n$ tangents from the Darboux Theorem \ref{th:darboux1}
form a Poncelet $n$-polygon inscribed in a conic $C$,  Darboux
proved that the curve $S$ of degree $n-1$ then completely
decomposes. More precisely, Darboux proved the following

\medskip

\begin{theorem} [Darboux, \cite {Dar1}]\label{th:darboux2}
 If a curve $S$ of degree $n-1$ is Poncelet-Darboux
related to a conic $K$ and if there is a conic $C$, a component of
$S$ which is $n$-Poncelet related to the conic $K$, then for
$n=2k+1$ the curve $S$ is completely decomposed on $k$ conics and
if $n=2k$ it is decomposed on $k-1$ conics and a line.
\end{theorem}

\medskip

In the case of decomposition of Poncelet-Darboux curve, the conic
components are parts of what we called the Poncelet-Darboux grids.
Further generalizations of Darboux theorems and Poncelet-Darboux
grids are obtained very recently in \cite{DR2}.

\medskip

Among other modern investigations in framework of the Darboux
theorems, we should mention here \cite{Tr}, \cite{Tr1}, \cite{NT}
and references therein, where  Poncelet-Darboux curves are related
to the study of stable bundles, instanton bundles and their
decomposition.

We reformulate the Darboux Theorem \ref{th:darboux1} following
\cite{Tr}:

\medskip

\begin{theorem} [Darboux]
Let $K$, $S$ be a nondegenerate conic and a curve of degree $n-1$
in the projective plane $\mathbb PW$ and let $\beta:W\mapsto W^*$
be a nondegenerate bilinear form. Assume that there are $n$ points
on $K$ such that the points of intersection of any two lines
associated to these points by $\beta$ belong to $S$. Then, $S$ is
Poncelet-Darboux $n-1$ related to $K$.
\end{theorem}

\medskip

\subsection{Some notations and notions}

\medskip

We will use the following notations.

By
\begin{equation}\label{eq:multiset}
\{x_1, x_2,\dots, x_n\}_M
\end{equation}
we will denote a {\it multiset}, meaning that number of
appearances of an item is important, but the order no; the notion
of divisor has synonymous meaning. The standard symmetric
functions of $n$ quantities $(a_1,\dots, a_n)$ will be denoted as
\begin{equation}\label{eq:symmetric}
\aligned
\sigma_1(a_1,\dots,a_n)&=\sum_{i=1}^na_i;\\
\sigma_2(a_1,\dots,a_n)&=\sum_{i<j}^na_ia_j;\\
\dots \\
\sigma_n(a_1,\dots,a_n)&=a_1a_2\cdots a_n;
\endaligned
\end{equation}
when one of the quantities, $a_k$, is omitted, the symmetric
functions of the $n-1$ rest quantities will be denoted as
\begin{equation}\label{eq:symmetric}
\aligned
\sigma_0^k(a_1,\dots,a_n)&=1 \quad k=1,\dots,n\\
\sigma_1^k(a_1,\dots,a_n)&=\sum_{i\ne k}^na_i  \quad k=1,\dots,n;\\
\sigma_2^k(a_1,\dots,a_n)&=\sum_{i<j, i,j\ne k}^na_ia_j  \quad k=1,\dots,n;\\
\dots \\
 \sigma_n^k(a_1,\dots,a_n)&=a_1a_2\cdots
a_{k-1}a_{k+1}\cdots a_n \quad k=1,\dots,n.
\endaligned
\end{equation}
\medskip
We will use also vector notation
\begin{equation}\label{eq:vector}
\overrightarrow\sigma_i(a_1,\dots , a_n)=(\sigma_i^1(a_1,\dots ,
a_n), \sigma_i^2(a_1,\dots , a_n), \dots, \sigma_i^n(a_1,\dots ,
a_n)),
\end{equation}
for $i=0, 1, \dots, n$;
\begin{equation}\label{eq:vector1}
\aligned
&\overrightarrow m=(m_1,\dots, m_n),\\
&\langle \overrightarrow m, \overrightarrow\sigma_i \rangle =
\sum_{k=1}^n m_k\sigma_i^k.
\endaligned
\end{equation}

\medskip
Let us recall some traditional notions: {\it the cyclic points}
are points on the infinite line with coordinates $\hat I=(1, i,
0)$ and $\hat J=(1,-i,0)$. A line is {\it isotropic} or {\it
minimal} if it is finite and passes through one of the cyclic
points. For a given curve, a point is {\it focal} if it is
intersection of two isotropic tangents to the curve with finite
points of contact with the curve. As an example, an ellipse has
four focal points: a pair of real foci and a pair of imaginary
foci.

\medskip

\section{Isofocal deformations}\label{sec:isofocal}

\medskip

\subsection{Definition of an integrable dynamical system}

\medskip

Let us start with a function of the form
\begin{equation}\label{eq:polynomial}
P^0(z)=(z-\alpha_1^0)^{m_1^0}
(z-\alpha_2^0)^{m_2^0}\cdots(z-\alpha_n^0)^{m_n^0},
\end{equation}
where all $\alpha_k^0$ are distinct, and consider its logarithmic
derivative  $F^0(z)=LP^0(z)=d[\log P^0(z)]/dz$:
\begin{equation}\label{eq:function}
F^0(z)=\frac{m_1^0}{z-\alpha_1^0}+\dots
+\frac{m_n^0}{z-\alpha_n^0}.
\end{equation}

Let us set
\begin{equation}\label{eq:function}
F^0(z)=\frac{m_1^0}{z-\alpha_1^0}+\dots
+\frac{m_n^0}{z-\alpha_n^0}=\frac{f(z)}{\phi(z)},
\end{equation}
where
\begin{equation}\label{eq:fphi}
\aligned \phi(z)&=(z-\alpha_1^0)(z-\alpha_2^0)\cdots
(z-\alpha_n^0)\\
f(z)&=\sum_{i=1}^n m_i^0\prod_{j\ne i}(z-\alpha_j^0).
\endaligned
\end{equation}
\medskip
One can easily see that
\begin{equation}
f(z)=B_nz^{n-1}+\dots +B_iz^{i-1}+\dots+B_1
\end{equation}
where, using notations from eq. (\ref{eq:vector}) and eq.
(\ref{eq:vector1}), we have
\begin{equation}\label{eq:integral1}
B_i=\langle \overrightarrow m^0, \overrightarrow
\sigma_{n-i}(\alpha_1^o,\dots, \alpha_n^0)\rangle.
\end{equation}
\medskip
We point out two particular cases.
\medskip
\begin{lemma}\label{lemma:derconst}
\begin{itemize}
\item[(a)] The function $f$ is equal to the derivative of $\phi$
if and only if all $m_i^0$ are equal to 1. \item[(b)] The function
$f$ is constant if and only if
\begin{equation}\label{eq:const}
\overrightarrow m^0\bot [\overrightarrow \sigma_0, \overrightarrow
\sigma_1, \dots, \overrightarrow \sigma_{n-2}].
\end{equation}
\end{itemize}
\end{lemma}

\medskip
Before proceeding with introduction of dynamics, we are going to
consider the simplest case as an example.

\medskip

\begin{example}
Let us consider the case $n=3$. According to the Marden Theorem
\ref{th:marden} for $n=3$, there exist a Marden curve $K$, which
is in this case,  {\it a conic}. Its   focal points $z_1$ and
$z_2$ satisfy  $f(z_i)=0$, under the condition $\deg f=2$. Again,
by the Marden Theorem \ref{th:marden}, the conic $K$ touches
line-segments $[\alpha_i^0, \alpha_j^0]$ in the ratio
$m_i^0:m_j^0$.

Now, since the lines $(\alpha_i^0, \alpha_j^0)$ are tangent to the
{\it conic} $K$, we may apply the Darboux Theorem
\ref{th:darboux1}. The triplet of the  conic $K$ and the two
polynomials $\phi$ and $f$ uniquely determines the
Poncelet-Darboux curve $\mathcal {PD}_K(\alpha_1^0, \alpha_2^0,
\alpha_3^0, m_1^0, m_2^0, m_3^0)=\mathcal {PD}_K(\phi, f)$. The
curve $\mathcal {PD}_K(\phi, f)$ is a {\it conic}.

Thus, the conics $\mathcal {PD}_K(\phi, f)$ and $K$ are 3-Poncelet
related.

According to the theorems of Poncelet and Darboux, there exists
another set of three points $\alpha_1^1, \alpha_2^1, \alpha_3^1$
which belong to the Poncelet-Darboux conic $\mathcal {PD}_K(\phi,
f)$ such that the lines $(\alpha_i^1, \alpha_j^1)$ are tangent to
the Marden conic $K$. The triangle $\alpha_1^1, \alpha_2^1,
\alpha_3^1$ is a Poncelet triangle with the caustic $K$ and the
boundary $\mathcal {PD}_K(\phi, f)$.

Now, we want to apply the Marden Theorem \ref{th:marden} with a
new polynomial
$$
\phi_1(z)=(z-\alpha_1^1)(z-\alpha_2^1)(z-\alpha_3^1)
$$
instead of $\phi$. In order to do that, we determine  {\it new
"masses"} $(m_1^1, m_2^1, m_3^1)$ such that the polynomial $f$
rests unchanged. More precisely, we calculate $(m_1^1, m_2^1,
m_3^1)$, up to a scalar factor, from the system of linear
equations:
\begin{equation}
\aligned B_3&=m_1^1+m_2^1+m_3^1\\
B_2&=m_1^1(\alpha_2^1+\alpha_3^1)+m_2^1(\alpha_1^1+\alpha_3^1)+m_3^1(\alpha_1^1+\alpha_2^1),\\
B_1&=m_1^1(\alpha_2^1\alpha_3^1)+m_2^1(\alpha_1^1\alpha_3^1)+m_3^1(\alpha_1^1\alpha_2^1)
\endaligned
\end{equation}
where the constants $B_1, B_2, B_3$ are determined from the
equation (\ref{eq:integral1})
$$
B_i=\langle \overrightarrow m^0, \overrightarrow
\sigma_{n-i}(\alpha_1^0,\alpha_2^0, \alpha_3^0)\rangle.
$$
Since $f$ is unchanged, the focal points of the Marden curve in
the new case are the same as focal points of $K$. We deduce now
that the Marden curve of the new case is equal to $K$, because
among confocal conics there is at most one inscribed in a
triangle.

By a new application of the Marden theorem, we finally get a new
information concerning Poncelet triangles. We are able to deduce
the ratio of tangency of a new triangle by the caustic $K$:
\begin{equation}\label{eq:newratio}
\aligned \alpha_3^1\beta_2^1:\beta_2^1\alpha_1^1&=m_3^1:m_1^1,\\
\alpha_2^1\beta_1^1:\beta_1^1\alpha_3^1&=m_2^1:m_3^1,\\
\alpha_1^1\beta_3^1:\beta_3^1\alpha_2^1&=m_1^1:m_2^1,
\endaligned
\end{equation}
where $\beta_i^1$ denote points of contact of the caustic and the
triangle.

\medskip

If the degree of the polynomial $f$ is equal to 1, then the conic
$K$ is a parabola. Thus we finish the Example with a nice
classical Lemma:

\medskip

\begin{lemma}[folklore]
If parabola touches three lines of a triangle $ABC$, then the
focus belongs to the circumscribed circle to the triangle $ABC$.
\end{lemma}
\medskip
If the polynomial $f$ is constant, then the conic $K$ is a circle.

\end{example}

\medskip

The last Example, although treatise the simplest case $n=3$ is
very instructive. In order to extract a new statement about ratios
of tangency of a new Poncelet triangle as it is formulated in eq.
(\ref{eq:newratio}), one needs several nontrivial steps of
alternative use of the Marden and the Poncelet- Darboux theorems.
In case $n>3$ we cannot follow the same lines, because the Marden
curve is not a conic any more and one cannot build up the Darboux
theorem straight on it. In order to link together the Marden and
the Darboux Theorems in general case we need to develop much more
subtle approach.

From the last Example, we learnt that it was fruitful idea to pass
from one Poncelet triangle to a new one by
\medskip
{\it a transformation which keeps the polynomial $f$ unchanged.}
\medskip

Following this principle, we introduce a new dynamics, which
depends on continuous "time" parameter $t$ and which has
quantities $(\alpha_1^o,\dots, \alpha_n^0)$ and $(m_1^0,\dots,
m_n^0)$ as {\it the initial data}. More precisely, we introduce
functions:
\medskip
\begin{equation}\label{eq:defdynam1}
\aligned \alpha_1(t)&=\alpha_1(t, \alpha_1^o,\dots, \alpha_n^0,
m_1^0,\dots, m_n^0),\\
\alpha_2(t)&=\alpha_2(t, \alpha_1^o,\dots, \alpha_n^0,
m_1^0,\dots, m_n^0),\\
\dots \\
\alpha_n(t)&=\alpha_n(t, \alpha_1^o,\dots, \alpha_n^0,
m_1^0,\dots, m_n^0),\\
m_1(t)&=m_1(t, \alpha_1^o,\dots, \alpha_n^0,
m_1^0,\dots, m_n^0),\\
m_2(t)&=m_2(t, \alpha_1^o,\dots, \alpha_n^0,
m_1^0,\dots, m_n^0),\\
\dots \\
m_n(t)&=m_n(t, \alpha_1^o,\dots, \alpha_n^0, m_1^0,\dots, m_n^0)
\endaligned
\end{equation}
in order to satisfy
\begin{equation}\label{eq:defdynam2}
F_t(z):=\frac{m_1(t)}{z-\alpha_1(t)}+\dots
+\frac{m_n(t)}{z-\alpha_n(t)}=\frac{f(z)}{\phi(z)+tf(z)},
\end{equation}
with the initial conditions
\begin{equation}\label{eq:defdynam3}
\aligned
m_i(o)&=m_i^0,\quad i=1,\dots, n,\\
\alpha_i(0)&=\alpha_i^0.
\endaligned
\end{equation}

\medskip
By the condition (\ref{eq:defdynam2}), the function $f$ keeps
unchanged during the evolution. This means that focal points, as
zeros of the polynomial $f$ are fixed during the evolution.Thus,
we will refer to such dynamics as {\it isofocal dynamics} or {\it
isofocal deformations}.
\medskip

The polynomial $f$ plays a role of an isospectral polynomial. From
the formula
\begin{equation}
f(z)=\sum_{i=1}^n m_i(t)\prod_{j\ne i}(z-\alpha_j(t)).
\end{equation}
we get the following
\medskip
\begin{proposition} The dynamics (\ref{eq:defdynam1},
\ref{eq:defdynam2}, \ref{eq:defdynam3}) has the coefficients of
the polynomial $f$ as first  integrals:
\begin{equation}\label{eq:integral}
B_i=\langle \overrightarrow m(t), \overrightarrow
\sigma_{n-i}(\alpha_1(t),\dots, \alpha_n(t))\rangle, \quad
i=1,\dots, n.
\end{equation}
\end{proposition}
\medskip
We will use terminology {\it positions in a moment $t$} for
$\alpha_1(t),\dots, \alpha_n(t))$ and for $(m_1(t),\dots, m_n(t))$
we will use {\it masses} although these masses will change during
the time and might be negative as well.

One of the first integrals is {\it the law of conservation of
masses}.

We will also use notation
\begin{equation}\label{eq:Phi}
\Phi_t(z):=\phi(z)+tf(z)=(z-\alpha_1(t))\cdots(z-\alpha_n(t)).
\end{equation}

\medskip

\subsection{Moser's trick and the Flashka coordinates}
\medskip

We consider the function (\ref{eq:defdynam1})
$$
F_t(z)=\frac{f(z)}{\Phi_t(z)}
$$
and we apply the Moser trick (see \cite{Mo}, \cite{Mo1}) to
develop it in a continued fraction of the following form.
\begin{equation}\label{eq:cfrac}
F_t(z)=\frac{1}{z-b_n-\frac{a_{n-1}^2}{z-b_{n-1}-\frac{a_{n-2}^2}{z-b_{n-2}-\dots
\frac{}{\dots -\frac{a_1^2}{z-b_1}}}}}
\end{equation}
The last formula (\ref{eq:cfrac}) gives us transformation from our
dynamical coordinates $(\alpha_1,\dots,\alpha_n, m_1,\dots, m_n)$
to the new coordinates $(a_1, a_2,\dots, a_{n-1}, b_1, b_2,\dots,
b_n)$. The last set of coordinates we will call {\it the Flashka
coordinates} (see \cite {Fl}, \cite{Mo},\cite{Mo1}).

\medskip

To construct the inverse transformation we consider the Flashka
Lax matrix $L_n$ for the $n$-point Toda chain (see \cite {Fl},
\cite{Mo},\cite{Mo1}):
\begin{equation}\label{eq:lax}
L_k=\left[\begin{array}{lllll}
b_1 & a_1 & & & \\
a_1 & b_2 & & & \\
    &     &\cdots & &\\
    &     &       & b_{k-1} & a_{k-1}\\
        &     &   & a_{k-1}& b_k
\end{array}\right]
\end{equation}
\medskip
Denote
$$
\delta_k=\det L_k.
$$
The following well-known difference relations take place
\begin{equation}\label{eq:difference}
\delta_k=(z-b_k)\delta_{k-1}-a_{k-1}^2\delta_{k-2},\quad
k=3,\dots,n.
\end{equation}
\medskip
The inverse transformation from the Flashka coordinates to the
initial dynamical coordinates is defined by the formula (see
\cite{Mo}):
\begin{equation}\label{eq:inverse}
F_t(z)=\frac{\delta_{n-1}}{\delta_n}.
\end{equation}
\medskip
From the last formula and from \ref{eq:defdynam2} we conclude
\medskip
\begin{lemma}\label{lemma:delta}
The time evolution of $\delta$ according to dynamics
\ref{eq:defdynam2} satisfies:
\begin{equation}\label{eq:difference1}
\aligned \delta_{n-1}(t)&=\delta_{n-1}(0)\\
\delta_n(t)&=\delta_n(0) + t\delta_{n-1}(0).
\endaligned
\end{equation}
\end{lemma}
\medskip
From the Lemma (\ref{lemma:delta}) one concludes that
\begin{itemize}
\item[$1^0$] {\it only $b_n$ among the Flashka coordinates depends
on $t$}; \item[$2^0$] {\it the coordinate $b_n$ depends on $t$
linearly.}
\end{itemize}
Thus we have the following
\medskip
\begin{theorem}\label{th:flashka}
The dynamical system (\ref{eq:defdynam2}) trivializes in the
Flashka coordinates, where it gets the form
\begin{equation}\label{eq:flashka}
\aligned \dot a_1=0 \quad \dot b_1&=0 \\
\dot a_2=0 \quad \dot b_2&=0 \\
\dots \\
\dot a_{n-1}=0 \quad \dot b_{n-1}&=0\\
 \dot b_n&=-1.
 \endaligned
\end{equation}
\end{theorem}
\medskip

\subsection{From Marden's curve to Poncelet-Darboux curve}
\label{subsection:mardenpondar}
\medskip

Following Marden, denote $\mathcal L_{\alpha_j^0}$ the line
equation of the point $\alpha_j^0=x_j+i y_j$:
$$
\mathcal L_{\alpha_j^0}=\lambda x_j +\mu y_j -1,
$$
the equation of all lines passing through the point $\alpha_j^0$.
The equation of the Marden curve (see \cite {Ma}, equation (4.10))
is deduced from the condition
\begin{equation}\label{eq:mardencurve}
\sum_{j=1}^n\frac{m_j^0}{\mathcal L_{\alpha_j^0}}=0.
\end{equation}
\medskip
We will denote the last curve $\mathcal M(\alpha_1^0, \alpha_2^0,
\dots,\alpha_n^0, m_1^0, m_2^0,\dots, m_n^0)$.
\medskip
Now we pass to projective plane. We see it as projective space of
quadratic polynomials: to a polynomial $P(z)=a(z-b)(z-c)$ we
associate the point $(-(c+b), cb, 1)$. The projection
\begin{equation}\label{eq:projection}
\pi: \mathbb {CP}^1\times \mathbb {CP}^1 \rightarrow \mathbb
{CP}^2
\end{equation}
is branched over the conic $K$ with the equation
\begin{equation}\label{eq:conic}
z_1^2=4z_0z_2.
\end{equation}
\medskip
In this presentation, to a linear polynomial $a(z-\alpha_j^0)$
corresponds the line $t_K(\alpha_j^0)$, tangent to the conic $K$.
This line is the set of all quadratic polynomials which have
polynomial $z-\alpha_j^0$ as a factor.

To develop further this connection, following Darboux (see
\cite{Dar1}) we introduce new system of coordinates. Given a plane
with standard coordinates $(z_0, z_1, z_2)$, we start from the
given conic $K$. It is given by the equation (\ref{eq:conic}) and
rationally parameterized by $(s^2, 2s, 1)$. The tangent line to
the conic $K$ through the point with the parameter $s_0$ is given
by the equation
$$
t_K(s_0): z_2s_0^2+z_1s_0+z_0=0.
$$
On the other hand, for a given point $P$ in the plane with
coordinates $P=(\hat z_0,\hat z_1,\hat z_2)$ there correspond two
solutions $\rho$ and $\rho_1$ of the quadratic in $s$ equation
\begin{equation}\label{eq:darboux}
\hat z_2s^2+\hat z_1s+\hat z_0=0.
\end{equation}
Each solution correspond to a tangent to the conic $K$ from the
point $P$. We will call the pair $(\rho, \rho_1)$ {\it the Darboux
coordinates} of the point $P$. One finds immediately
\begin{equation}\label{eq:viet}
\frac{\hat z_0}{\rho \rho_1}=-\frac{\hat z_1}{\rho+\rho_1}=\hat
z_2.
\end{equation}
\medskip
The line $t_K(\alpha_j^0)$ which corresponds to the point
$\alpha_j^0$ and to the linear polynomial $a(z-\alpha_j^0)$ has
the following presentation in the Darboux coordinates:
$$
t_K(\alpha_j^0)(\rho,
\rho_1)=(\rho-\alpha_j^0)(\rho_1-\alpha_j^0).
$$
The Poncelet-Darboux curve is done by the equation
\begin{equation}\label{eq:darbouxcurve}
\sum_{j=1}^n \frac{m_j^0}{t_K(\alpha_j^0)}=0.
\end{equation}
In Darboux coordinates it may be rewritten in the form
\begin{equation}\label{eq:darbouxcurve2}
\frac{f(\rho)}{\phi(\rho)}=\frac{f(\rho_1)}{\phi(\rho_1)}.
\end{equation}
\medskip
The curve defined by last equations we will denote as $\mathcal
{PD}_K(\phi, f)$. From the equations (\ref{eq:mardencurve}) and
(\ref{eq:darbouxcurve}) and whole previous considerations we get
the following theorem.
\medskip
\begin{theorem}\label{th:mardendarboux} There is a
birational morphism defined by the equations (\ref{eq:fphi}) and
(\ref{eq:integral}) between the data of Marden curves
$(\alpha_1^0, \alpha_2^0,\dots, \alpha_n^0, m_1^0, m_2^0,\dots,
m_n^0)$ and the data $(\phi, f)$ of Poncelet-Darboux curves
associated with the conic $K$. There is a birational morphism
between dual of projective closure of a Marden curve  $\mathcal
M(\alpha_1^0, \alpha_2^0,\dots, \alpha_n^0, m_1^0, m_2^0,\dots,
m_n^0)$ and the Poncelet-Darboux curve $\mathcal {PD}_K(\phi, f)$
with corresponding data.
\end{theorem}
\medskip

\medskip
\subsection{Discriminant and gauge equivalence}\label{sec:discrim}
\medskip

The morphism from the previous Theorem (\ref{th:mardendarboux})
fails to be one to one for those $\alpha_1^0, \alpha_2^0,\dots,
\alpha_n^0$ for which the system of linear equations in  $(m_1^0,
m_2^0,\dots, m_n^0)$ given by the equations (\ref{eq:integral})
has determinant $D(\alpha_1^0, \alpha_2^0,\dots, \alpha_n^0)$
equal zero. The condition
$$
D(\alpha_1^0, \alpha_2^0,\dots, \alpha_n^0)=0
$$
is equivalent to $\alpha_i^0=\alpha_j^0$ for some $i\ne j$.

We will refer to configurations $(\alpha_1(t), \alpha_2(t),\dots,
\alpha_n(t), m_1(t), m_2(t),\dots, m_n(t))$ such that
$\alpha_i(t)=\alpha_j(t)$ for some $i\ne j$ as {\it points of
collision}. Such a moment $t$ we will call {\it moment of
collision}. If there is no $k\ne i, k\ne j$ such that in addition
$\alpha_i(t)=\alpha_j(t)=\alpha_k(t)$ the point of collision is
{\it simple}. The system is {\it simple} if it has only simple
collision points.
\newline
{\it We will assume that system passes smoothly through a
collision point  in the phase space.}
\medskip

\section{Complete decomposition of Poncelet-Darboux
curves, $n$-volutions and collisions}\label{sec:decomposition}
\medskip

\subsection{$n$-volutions}
\medskip
By an {\it $n$-volution} in a set $V$ we will assume a family of
multisets $\mathcal A_n$, a subset of the $n$-th symmetric product
$Sym_n V$ such that there is a unique function $f$:
$$
f:V\rightarrow \mathcal A_n,
$$
such that $v\in f(v)$. In some classical terminology is used
notion of cyclic-symmetric correspondences.

If $\alpha\in \mathcal A_n$ is a multiset such that its
cardinality is less than $n$ we will say that it is {\it a
collision point} of the $n$-volution. In  other words
$\alpha=\{\alpha_1,\dots,\alpha_n\}_m$ is a collision point if
there exist $\alpha_1,\dots,\alpha_k$, $k<n$ and natural numbers
$c_1,\dots, c_k$ such that
$$
\alpha=\{\alpha_1,\dots,\alpha_n\}_m=c_1\alpha_1+\dots
c_k\alpha_k.
$$
The basic examples of $n$-volutions are involutions, which
correspond to $n=2$. Then, the notion of collision point coincides
with the notion of fixed point. A nice case are involutions of a
conic. By the Fregier theorem, we know that for every involution
on a conic, there exists a point, {\it the Fregier point} such
that the involution is cut from the conic by lines from the pencil
determined by the Fregier point.

We pass now to the case of $\mathbb {CP}^1$. By Luroth Theorem we
know that every $n$-volution is determined by a pair of
polynomials $p, q$ of degree $n$.

Consider the pencil $p_t(z)=tp(z)+q(z)$ and the roots
$a_1(t),\dots,a_n(t)$. There is a one-parameter family of
$n$-tuples. From the following system
$$
\aligned tp(a_1(t))+q(a_1(t))&=0\\
tp(a_2(t))+q(a_2(t))&=0
\endaligned
$$
we get
$$
(a_1(t)-a_2(t))r(a_1(t),a_2(t))=0.
$$
Here $r$ is of degree $n-1$ in $a_1$ and symmetric. When $t$
varies, $r(a_1(t), a_2(t))$ defines a curve.

The question is how to describe all such $r$. Following the lines
of Section (\ref{subsection:mardenpondar}) we pass to $\mathbb
{CP}^2$ and correspond to a $n$-tuple of $n$ linear factors
$(z-a_1),\dots, (z-a_n)$ a polygon of $n$-sides circumscribed
about the conic $K$. Thus, we have the following
\medskip

\begin{proposition} An $n$-volution defined by polynomials $p, q$
of degree $n$ is associated to a curve from $Pon-Dar_{n-1}(K)$
given by the equation
$$
\det A\circ K(z)=0.
$$
The $(n-1)\times (n+1)$   matrix $A$ annihilates the pencil
generated by $p, q$ and the matrix $K(z)$ is induced by the conic
$K$ and it has the form
$$
K(z)=\left(\begin{array}{cccccc} z_2 & 0 & &\dots &\dots & 0\\
-z_1 & z_2 & & \dots &\dots & 0\\
z_0 & -z_2 & & \dots & \dots & 0\\
0 & z_0  & & \dots & \dots & \\
  &     & &\dots & & z_2\\
  &     & &\dots & & -z_1\\
  &     & &\dots & & z_0
\end{array}\right)
$$
\end{proposition}
\medskip

The last equation follows from \cite{Tr}. A linear system $L=L(p,
q)$ has a {\it base point} if polynomials $p$ and $q$ have a
common root. From \cite{Tr} we have that the base points of the
linear system correspond to the components of the Poncelet-Darboux
curve which are tangent lines to the conic $K$.  From \cite {Tr}
we also have the following description of the collision points.
\medskip
\begin{proposition}\label{prop:traut}
Let $C$ be a Poncelet-Darboux curve with respect to the conic $K$
such that the corresponding linear system $L=L(p, q)$ doesn't have
base points. Then for a given point $P=(s_0, t_0)\in K$ and an
integer $k\ge 0$ the following are equivalent:
\begin{itemize}
\item[(A)] the intersection multiplicity of $K$ and $C$ at $P$ is
equal to $k$; \item[(B)] $(s_0, t_0)$ is a zero of order $k$ of a
unique polynomial $h\in L(p,q)$.
\end{itemize}
\end{proposition}
\medskip
\begin{example}
Now, we want to calculate the number of {\it generalized Fregier
points} in the case of pencils of curves of degree $k$. The pencil
of curves of degree $k$ has a base set of $\binom {k+2}{2}-2$
points. A curve of degree $k$ intersects the conic $K$ in $n=2k$
points. But, between a polynomial of degree $n=2k$ and a curve of
degree $k$ which determine the same set of $2k$ points on the
conic $K$, there is a correspondence which is not bijective, due
to the relation
$$4z_0z_1=z_2^2$$
which holds on the conic $K$ by definition. Thus, one needs to fix
additionally $\binom {k}{2}$ points outside the conic to get
uniqueness. As a result, the number of generalized Fregier points
for an $n$-volution is
$$\binom {k+2}{2}-2-\binom {k}{2}=2k-1.$$
For the case $n=2k-1$ one needs also to fix a point on the conic
$K$, since the number of intersections of a curve and a conic is
always even. In this case the number of generalized Fregier points
is $2k-2$.

In both cases, with $n$ even or odd the total number of
generalized Fregier points is equal $n-1$. This number coincides
with the number of focal points in the Marden Theorem. Is there a
natural bijection between these two sets?
\end{example}

\medskip
\subsection{Conic component and complete decomposition of
Poncelet-Darboux curves}

\medskip
From the last Proposition (\ref{prop:traut}) we see that the
intersection of a Poncelet-Darboux curve $C$ and the conic $K$ is
transversal if and only if corresponding system is simple in
terminology of the Section (\ref{sec:discrim}). If a
Poncelet-Darboux curve $C$ satisfies any of two equivalent
conditions of the Proposition (\ref{prop:traut}), we will say that
it is {\it transversal}. Thus, the linear system which corresponds
to a transversal Poncelet-Darboux curve contains polynomials with
at most double roots.

From now on we will consider transversal Poncelet-Darboux curves.
Let us consider first the simplest case of curves of degree 2.
\medskip
\begin{example}\label{ex:conicponceletdarboux}
Let $C$ be degree two Poncelet-Darboux curve with transversal
intersection with the conic $K$. Denote four intersection points
$C\cap K=\{x_1, x_2, x_3, x_4\}$. Denote four common bitangents
$t_i$, $i=1,\dots, 4$ and denote points of contact of the tangent
$t_i$ with the conic $C$ as $y_i$ and the point of contact with
the conic $K$ as $a_i$.

The conic $C$ is 3-Poncelet related to the conic $K$. By the
Poncelet theorem (\ref{th:poncelet}) this means that there is a
triangle inscribed in $C$ and circumscribed about $K$ with
arbitrary point of $C$ taken as a vertex.

Choose any of the points $y_i$, say $y_1$ as a vertex. Then $y_1$
is a double vertex of a Poncelet triangle. The third vertex of the
Poncelet triangle is one of the points $x_i$ and the Poncelet
triangle is $y_1y_1x_i$. Thus we get
\medskip
\begin{lemma}\label{lemma:trhreeponcelet}
With use of previous notations, for any point $y_i$ there is a
point $x_j$ such that the line $y_ix_j$ is tangent to the conic
$K$ at the point $x_j$. The triplet $(y_iy_ix_j)$ forms a Poncelet
triangle.
\end{lemma}
\medskip
\end{example}
\medskip

Coming back to the general case of transversal Poncelet-Darboux
curves of degree $n-1$, let us study the case when a curve $S$ has
a conic component $C$ which is $n$-Poncelet related to the conic
$K$. According to the Darboux theorem (\ref{th:darboux2}) then the
curve $S$ is decomposed as a product of $k$ conics if $n=2k+1$ or
as a product of $k-1$ conics and a line if $n=2k$.

Any conic $C_i$ defines a symmetric $2-2$ correspondence of the
Euler-Chasles type $\Phi_i$ such that a point $P$ belongs to the
conic $C_i$ if and only if $\Phi_i(\rho, \rho_1)=0$ where $(\rho,
\rho_1)$ are the Darboux coordinates of the point $P$. Darboux
proved the Theorem (\ref{th:darboux2}) using these correspondences
together with the opposite statement  that a symmetric
correspondence determines a conic from the pencil. We know that
correspondences commute if the conics belong to the same pencil
with the conic $K$. Moreover, we have
\medskip
\begin{lemma}\label{lemma:transversal}
Suppose the lines $t_1, t_2,\dots, t_n$ are given such that
$\{P_i\}:=t_i\cap t_{i+1}\in C_i$ and $\{P_n\}:=t_n\cap t_1\in
C_n$, where conics $C_i$ belong to a confocal pencil together with
the conic $K$. Denote related Euler-Chasles correspondences as
$\Phi_i$. If multisets are equal
$$
\{\,C_2, C_3, \dots C_k\}_M=\{C_{n-1}, C_{n-2},\dots,
C_{n-k+1}\}_M
$$
then there exist a conic $C_0$ from the confocal family which
contains the point $P_1$ and the point $P_{k+1, n-k}:=t_{k+1}\cap
t_{n-k}$.
\end{lemma}
\medskip
\begin{proof}
The Lemma follows from the fact that in the plane there are two
conics from a confocal family which contain a point $P_1$. In our
notation one conic is $C_1$, denote the other one as $C_0$ with
the Euler-Chasles correspondence $\Phi_0$. We have:
\begin{equation}
\aligned \Phi_{k}\circ\dots \circ \Phi_2\circ
\Phi_1(t_1)&=t_{k+1}\\
\Phi_{k}\circ\dots \circ \Phi_2\circ \Phi_0(t_1)&=t_{k+1}.
\endaligned
\end{equation}
Using commuting property, we get from the last equation
\begin{equation}\label{eq:phi0}
\Phi_0\circ\Phi_{k}\circ\dots \circ \Phi_2(t_1)=t_{k+1}.
\end{equation}
Consider the intersection  of the lines
$$\Phi_{n-k+1}\circ\dots\Phi_{n-2}\circ\Phi_{n-1}(t_1)$$
and $t_{k+1}$. Denote by $\hat \Phi$ the Euler-Chasles
correspondence associated with the conic $C$ such that
$$
\Phi_{n-k+1}\circ\dots\Phi_{n-2}\circ \Phi_{n-1}(t_1)\cap
t_{k+1}\in C.
$$
Thus we have
$$
\hat \Phi\circ \Phi_{n-k+1}\circ\dots\Phi_{n-2}\circ
\Phi_{n-1}(t_1)= t_{k+1}.
$$
From the assumption of the Lemma, the last equation and from the
equation (\ref{eq:phi0}) it follows that $C=C_0$, giving the proof
of the Lemma.
\end{proof}
\medskip
There is an important special case of previous Lemma, when all the
conics are equal: $C_2=C_3=\dots
=C_k=C_{n-1}=C_{n-2}=\dots=C_{n-k+1}$. In the case of Poncelet
$n$-tangles even more is true: $C_i=C_j$ for any $i, j=1,\dots,
n$.
\medskip
The union of $k$ conics if $n=2k+1$ or $k-1$ conic and the line if
$n=2k$ together with conics from the previous Lemma
(\ref{lemma:transversal}) form {\it the complete projective
Poncelet-Darboux grid}. For the study of Poncelet-Darboux grids
see \cite{DR2} and references therein. Notice that transversal
conics from the Lemma (\ref{lemma:transversal}) doesn't form a
decomposition of a Poncelet-Darboux curve.
\medskip

\begin{example}
Suppose that all vertices of a triangle lie on a conic $C_1$, one
of the sides touches a conic $C_2$ and the other two sides touch
the conics $tC_1+C_2=0$ and $sC_1+C_2=0$ respectively. If the
three points of contact are not collinear, then
$$
(I_3-I_1ts)^2-4I_4(I_2+I_1(s+t))=0,
$$
where $I_1, I_2, I_3, I_4$ denote invariants of the pair of conics
$(C_1, C_2)$.
\end{example}

\medskip

\begin{example}
Let $P_1, P_2,\dots, P_n,\dots$ be points on a conic $C_1$ such
that there exists a conic $C_2$ to which all sides $P_iP_{i+1}$
are tangent. Assume $P_{k+2}\ne P_k$. Then that lines
$A_1A_{k+1}$, $A_2A_{k+2}$,\dots $A_iA_{i+k}$ touch the conic
$t_kC_1+C_2=0$ where
$$t_2=0, t_3=\frac{I_3^2-4I_2I_4}{4I_1I_4}$$
and
$$t_{k+1}=\frac{(I_3^2-4I_2I_3)-4I_1I_4t_k}{I_1^2t_k^2t_{k-1}}.$$
\end{example}

\medskip

The condition {\it for two conics $K$ and $C$} to be $n$-Poncelet
related has been derived by Cayley. Recent account of the subject
can be found in \cite {DR}, \cite {DR2}.

Here, we want to present a condition {\it for two polynomials}
$\phi, f$ of degree $n$ in order that corresponding
Poncelet-Darboux curve $S=\mathcal{PD}_K(\phi, f)\in
Pon-Dar_{n-1}(K)$ has a conic component $n$-Poncelet related to
the conic $K$.
\medskip

\begin{theorem}\label{th:condition}
Let a transversal Poncelet-Darboux curve:
$$
S=\mathcal{PD}_K(\phi,f)\in Pon-Dar_{n-1}(K)
$$
be given, where polynomials $f$ and $\phi$ are without common
zeros.
\begin{itemize}
\item[(i)] For $n=2k+1$, curve $S$ is completely decomposed to $k$
conics only if there exist four values $t_1, t_2, t_3, t_4$ such
that
\begin{equation}\label{eq:decomposition1}
\deg\GCD\left(\Phi_{t_i}(z),\frac{d}{dz}\Phi_{t_i}(z)\right)=k,
\end{equation}
for $i=1, 2, 3, 4$;

\item[(ii)] for $n=2k$, curve $S$ is completely decomposed to
$k-1$ conics and a line  only if there exist four values $t_1,
t_2, t_3, t_4$ such that:
\begin{equation}\label{eq:decomposition2a}
\deg\GCD\left(\Phi_{t_i}(z),\frac{d}{dz}\Phi_{t_i}(z)\right)=k,
\end{equation}
for $k=1,2$, and
\begin{equation}\label{eq:decomposition2b}
\deg\GCD\left(\Phi_{t_i}(z),\frac{d}{dz}\Phi_{t_i}(z)\right)=k-1,
\end{equation}
for $k=3,4$.
\end{itemize}
Here, we denoted $\Phi_{t_i}(z):=\phi(z)+t_if(z)$, while $\GCD$
stands for the greatest common divisor of polynomials.
\end{theorem}
\medskip
\begin{proof}
Suppose that Poncelet-Darboux curve $S$ of degree $n-1$ completely
decomposes. Then, denote by $C$ its conic component which is $n$
Poncelet related to the conic $K$. By assumption of
transversality, conics $C$ and $K$ intersect in four points. As in
Example (\ref{ex:conicponceletdarboux}) denote four intersection
points $C\cap K=\{x_1, x_2, x_3, x_4\}$. Denote four common
bitangents $t_i$, $i=1,\dots, 4$ and denote points of contact of
the tangent $t_i$ with the conic $C$ as $y_i$ and the point of
contact with the conic $K$ as $a_i$.

The conic $C$ is $n$-Poncelet related to the conic $K$ and
according to  the Poncelet theorem (\ref{th:poncelet}) there is a
polygon of  $n$ sides inscribed in $C$ and circumscribed about $K$
with arbitrary point of $C$ taken as a vertex.

Suppose $n$ is odd: $n=2k+1$. Choose any of the points $y_i$, say
$y_1$ as a vertex. Then $y_1$ is a double vertex of a Poncelet
$2k+1$-polygon. Moreover, next $k-1$ vertices $c_1, \dots,
c_{k-1}$ are also double vertices.  The last vertex of the
Poncelet $n$-tangle is one of the points $x_i$ and the Poncelet
$n$-tangle is $c_{k-1}\dots c_1y_1y_1c_1\dots c_{k-1}x_i$.

Suppose now that $n$ is even: $n=2k$. In this case there are two
pairs of distinguished Poncelet $n$-tangles. Two of them are of
the form $c_{k-1}\dots c_1y_jy_jc_1\dots c_{k-1}y_iy_i$, for
example for $(i, j)= (1,2)$ and for $(i, j)= (3, 4)$. Each of them
connect a pair of common tangents and it has all other vertices as
double. Another pair of distinguished Poncelet $n$-tangles connect
pair of intersection points. For example the first one connects
$x_1$ and $x_2$ while the second connects $x_3$ with $x_4$. All
other their vertices are double. Thus these two Poncelet
$n$-tangles  are of the form $d_{k-1}\dots d_1x_1d_1\dots
d_{k-1}x_2$ and $e_{k-1}\dots e_1x_3e_1\dots e_{k-1}x_4$.

\medskip
\begin{lemma} If a transversal Poncelet-Darboux curve $S=\mathcal{PD}_K(\phi,
f)\in Pon-Dar_{n-1}(K)$ has a conic component $C$ which is
$n$-Poncelet related to the conic $K$ then:
\begin{itemize}
\item[(i)] for $n=2k+1$  there exist four values $t_1, t_2, t_3,
t_4$ and four polynomials $Q_1, Q_2, Q_3, Q_4$ such that
\begin{equation}\label{eq:decomp1}
\aligned \Phi_{t_1}(z):=\phi(z)+t_1f(z)&=(z-a_1)Q_1^2(z)\\
\Phi_{t_2}(z):=\phi(z)+t_2f(z)&=(z-a_2)Q_2^2(z)\\
\Phi_{t_3}(z):=\phi(z)+t_3f(z)&=(z-a_3)Q_3^2(z)\\
\Phi_{t_4}(z):=\phi(z)+t_4f(z)&=(z-a_4)Q_4^2(z);
\endaligned
\end{equation}
\item[(ii)] for $n=2k$  there exist four values $t_1, t_2, t_3,
t_4$ and four polynomials $Q_1, Q_2, Q_3, Q_4$ such that
\begin{equation}\label{eq:decomp2}
\aligned \Phi_{t_1}(z):=\phi(z)+t_1f(z)&=(z-a_1)(z-a_2)Q_1^2(z)\\
\Phi_{t_2}(z):=\phi(z)+t_2f(z)&=(z-a_3)(z-a_4)Q_2^2(z)\\
\Phi_{t_3}(z):=\phi(z)+t_3f(z)&=Q_3^2(z)\\
\Phi_{t_4}(z):=\phi(z)+t_4f(z)&=Q_4^2(z).
\endaligned
\end{equation}
\end{itemize}
\end{lemma}
\medskip
From the conditions (\ref{eq:decomp1}) and (\ref{eq:decomp2})
immediately follow conditions (\ref{eq:decomposition1}) and
(\ref{eq:decomposition2a}) together with
(\ref{eq:decomposition2b}) of the Theorem respectively. This
proves  the Theorem.

\end{proof}

For the opposite direction, observe that from the conditions
(\ref{eq:decomposition1}) and (\ref{eq:decomposition2a}) together
with (\ref{eq:decomposition2b}) of the Theorem follow the
conditions of the Lemma (\ref{eq:decomp1}) and (\ref{eq:decomp2})
by use of the transversality condition. From the transversality it
follows that all multiple zeros of the polynomials in the pencil
generated by $f$ and $\phi$ are of the second degree. Now, from
the conditions (\ref{eq:decomp1}) and (\ref{eq:decomp2}) one can
easily prove the following
\medskip

\begin{lemma}
Suppose the conditions (\ref{eq:decomp1}) and (\ref{eq:decomp2})
are satisfied. Denote by $\Gamma_1$ and $\Gamma_2$ the following
elliptic curves:
\begin{equation}
\aligned \Gamma_1:& y^2=(z-a_1)(z-a_2)(z-a_3)(z-a_4)\\
\Gamma_2:& Y^2=(X+t_1)(X+t_2)(X+t_3)(X+t_4).
\endaligned
\end{equation}
Then, there is a $n:1$  morphism
$$
h:\Gamma_1 \rightarrow \Gamma_2,\quad h(z,y)=(X,Y),
$$
where
$$
X=\frac{\phi(z)}{f(z)},\quad
Y=y\frac{Q_1(z)Q_2(z)Q_3(z)Q_4(z)}{f^2(z)}.
$$
\end{lemma}
\medskip
From the last Lemma we see that question of decomposition of the
Poncelet-Darboux curve  defined by polynomials $f$ and $\phi$
corresponds to study of an unramified covering of degree $n$ of
the elliptic curve $\Gamma_1$ over the elliptic curve $\Gamma_2$.
Such a covering is realized as factorization of the first curve by
its finite subgroup. We can say even more about the above
covering.
\medskip

\begin{lemma}\label{lemma:covering}
\begin{itemize}
\item[(a)] There is a constant $N$ such that
$$
Q_1(z)Q_2(z)Q_3(z)Q_4(z)=N\left(\frac{d}{dz}\phi(z)f(z)-\phi(z)\frac{d}{dz}f\right).
$$
\item[(b)] There is a relation between holomorphic differentials
of elliptic curves $\Gamma_1$ and $\Gamma_2$:
$$
\frac{dz}{y} = N \frac {dX}{Y}.
$$
\end{itemize}
\end{lemma}
\medskip

The proof follows by straightforward calculations.

Denote by $\Lambda_1$ and $\Lambda_2$ lattices which correspond to
the elliptic curves $\Gamma_1$ and $\Gamma_2$ respectively. Denote
by $u$ a parameter on $\Gamma_1$. Then the covering defines
correspondence between functions $z(u|\Lambda_1)\mapsto
X(u/N|\Lambda_2)$. Thus we get
\medskip
\begin{lemma}\label{lemma:covering2}
For complete decomposition of the Poncelet-Darboux curve is
necessary that
$$
N = \frac{a}{n},
$$
and
$$ a\Lambda_2 \subset \Lambda_1.$$
\end{lemma}
\medskip
Then, if all above conditions are satisfied, one gets a Poncelet
trajectory by choosing $P\in \Gamma_1$ such that $nP\in
a\Lambda_2\subset \Lambda_1$. The points of the trajectory
correspond to parameters $u_j = u_0 + j\eta$, $j=0,\dots, n-1$,
where $\eta$ corresponds to the point $P$. Denote by $C_P$ a conic
which corresponds to $P$ in the pencil of conics generated with
the conic $K$ and its four tangents at the points $a_1, a_2, a_3,
a_4$. Then the conic $C_P$ is $n$-Poncelet related to the conic
$K$.
\medskip

The systematic study of elliptic coverings was established in
\cite{Jac2}. We pass now to Jacobi's notation. By applying
rational-linear transformations, we come to the canonical form of
elliptic curves $\Gamma_2$ and $\Gamma_1$ and relation between
differentials
$$
\frac {N\,dy}{\sqrt{(1-y^2)(1-\lambda
y^2)}}=\frac{dx}{\sqrt{(1-x^2)(1-k\cdot x^2)}}.
$$
The constants $k$ and $\lambda$ are the modules of the elliptic
curves $\Gamma_1$ and $\Gamma_2$.

Following Jacobi, for given  $n$ odd and given module $k$, we have
explicit formulae for the transformations.
\medskip
\begin{theorem}\label{th:jacobitransformations}
For given $n$ odd and for
$$
\omega = \frac{m K + m' K'}{n},
$$
where integers $m$ and $m'$  have no common divisors which divide
$n$, the transformation is defined by
$$
\aligned
f(z)&=\frac{x}{N}\prod_{r=1}^{(n-1)/2}\left(1-\frac{x^2}{\sn^24r\omega}\right)\\
\phi(z) &=\prod_{r=1}^{(n-1)/2}(1-x^2 \cdot\sn^24r\omega)\\
N &=(-1)^{(n-1)/2}\prod_{r=1}^{(n-1)/2}\left(1-\frac{\sn(K-4r\omega)}{\sn^2(4r\omega)}\right)\\
\lambda &= \prod_{r=1}^{(n-1)/2}\left(\sn^4(K-4r\omega)\right).
\endaligned
$$
The transformation corresponds to an $n$- Poncelet trajectory,
where
$$
x_i=\sn(u+4(i+1)\omega, k)\quad i=0,\dots, n-1, \quad y=\sn(u/N,
\lambda).
$$
\end{theorem}
\medskip
We can make another view on the situation. Let us consider three
arithmetic functions. Suppose an natural number $n\in\mathbb N$ be
given by its prime decomposition:
$$
n = p_1^{k_1}p_2^{k_2}\cdots p_r^{k_r},
$$
where $p_i$ are different prime numbers. Following \cite {BM}, we
define a function $t(n)$, {\it the number of primitive $n$-torsion
points on an elliptic curve}:
$$
t(n):=(p_1^2-1)p_1^{2(k_1-1)}(p_2^2-1)p_2^{2(k_2-1)}\cdots
(p_r^2-1)p_r^{2(k_r-1)}.
$$
As an example, for $n=p$ a prime number, $t(p)=p^2-1$. The
function $t$ is a multiplicative arithmetic function. It
represents the number  of Poncelet polygons in total up to the
porism, with fixed caustic and confocal pencil of conics.

As the second arithmetic function, we introduce a function
$\sigma'(n)$ as a multiplicative function which is for $n$ odd
equal to
$$
\sigma'(n)=n \left (1 + \frac{1}{p_1}\right)\left (1
 + \frac{1}{p_2}\right)\cdots \left (1 + \frac{1}{p_r}\right),\quad
 n\,\,
odd.
$$
For $n=2^k$ we define
$$
\sigma'(2^k):=2^{k+1}-1.
$$
For example, for $n=p$ a prime number, we have
$$
\sigma'(p) = p + 1 = \sigma (p),
$$
the $\sigma' $ function in this case is equal to the $\sigma$
function, the sum of divisors of $p$. The function $\sigma' (n)$
counts the number of degree $n$ elliptic coverings of the above
form assuming the module $k$ being fixed. The number of
transformations listed in Theorem (\ref{th:jacobitransformations})
for given odd number $n$ is equal to $\sigma' (n)$.

The third arithmetic function we are going to consider is  well
known Euler function $\varphi (n)$ counting the numbers smaller
than $n$  relatively prime to $n$:
$$
\varphi (n)= n \left (1 - \frac{1}{p_1}\right)\left (1
 - \frac{1}{p_2}\right)\cdots \left (1 - \frac{1}{p_r}\right).
$$
\medskip
\begin{proposition} For $n$ odd the identity holds:
$$
t(n) = \sigma' (n) \cdot \varphi (n).
$$
\end{proposition}
\medskip
\begin{proposition} For $n=2^k$, $k\ge 1$, the inequality
holds:
$$
t(n)\le \sigma' (n) \cdot \varphi (n).
$$
The equality holds only for $k=1$.
\end{proposition}
\medskip
\begin{theorem}\label{th:cyclic} Let $m$ be an arbitrary odd number. For $n=2^k \cdot m$
where $k=0, 1$ all above elliptic coverings are cyclic.

For $n=2^k \cdot m$ and every $k>1$, there are elliptic coverings
of the above form which are not cyclic.
\end{theorem}
\medskip
Thus, we come to the converse of the Theorem (\ref{th:condition}).
\medskip
\begin{theorem}\label{th:convereseconditions}
Let $m$ be an arbitrary odd number. For $n=2^k \cdot m$ where
$k=0, 1$ the conditions of the Theorem (\ref{th:condition}) are
sufficient as well.
\end{theorem}
\medskip
For $n=2^k \cdot m$, where $k>1$, one needs to do careful analysis
to distinguish those coverings which are cyclic. The general
description of transformations of an even order $n$ is given by
$$
\aligned \phi(z)&=\frac{1}{2}((1+z)(1+kz)T^2+(1-z)(1-kz)T'^2)\\
f(z)&=\frac{1}{2}((1+z)(1+kz)T^2-(1-z)(1-kz)T'^2),
\endaligned
$$
where
$$
T(z) = P(z) + zQ(z), \quad T'(z) = P(z) - zQ(z).
$$
Here $P, Q$ are even polynomials, in other words, polynomials  in
$z^2$.

Suppose a transformation of degree $n$ consists of $k$ cycles of
length $l$, where $n=kl$:
$$
\begin{array}{cccc}
\sn u_0 & \sn (u_0 + \omega)  & \dots & \sn (u_0 + (l-1)\omega)\\
\sn u_1 & \sn (u_1 + \omega)  & \dots & \sn (u_1 + (l-1)\omega)\\
\dots & \dots & \dots & \dots \\
\sn u_{k-1} & \sn (u_{k-1} + \omega)  & \dots & \sn (u_{k-1} +
(l-1)\omega)
\end{array}
$$
where
$$u_i = u_{i-1} + \omega_i, \quad i=1, \dots k-1,$$
and
$$\omega_k=\omega - \sum_{i=1}^{k-1}\omega_i.$$
Denote by $C_i$ the conic which corresponds to $\omega_i$ in the
pencil of conics defined by $K$ and the four tangents at $a_j$,
$j=1,\dots, 4$. Geometric interpretation of the last
transformation is realized through the Complete Poncelet Theorem:
corresponding Poncelet polygon consists of $n=kl$ tangents to the
conic $K$ with vertices
$$P_1P_2\dots P_kP_{k+1}\dots P_{2k}\dots
P_{(l-1)k+1}\dots P_{lk}$$ where $P_1\in C_1$, $P_2\in C_2$, and
more generally
$$P_s\in C_r \Leftrightarrow s\equiv r \,(\mod k).$$

The same interpretation may be done for cyclic transformations of
composite degree $n=kl$ as well. In the case of cyclic
transformations of composite order, in  this way we get another
geometric realization beside the one connected with Poncelet
theorem from the beginning.

\medskip
\subsection{Conic choreography of the "$\binom {n}{2}$-body" problem}
\label{subsec:choreography}
\medskip

In this Section we want  to get effective description of
polynomials $f$ and $\phi$ which satisfy the  Theorem
(\ref{th:condition}). We use parametrizations of the
transformations obtained by classics algebraically. For small odd
numbers, we give complete description of initial data of the
isofocal transformations which correspond to the cases of complete
decomposition of the Poncelet-Darboux curves.
\medskip

Assume $n$ is an odd number. Then all transformations are given by
the formulae
$$
\aligned \phi(x)&= P^2+2PQ+Q^2x^2\\
f(x)&= x(P^2+2PQ+Q^2x^2)
\endaligned
$$
where $P$ and $Q$ are polynomials in $x^2$:
$$
\aligned P(x)&=\alpha + \gamma x^2 + \epsilon x^4 +\dots\\
Q(x)&= \beta + \delta x^2 + \zeta x^4+\dots,
\endaligned
$$
of degree $p$ in $x^2$ both  if $n=4p+3$ and of degree $p$ and
$p-1$ respectively if $n=4p+1$. Then
$$
\frac {1}{N}= 1 + \frac{2\beta}{\alpha}.
$$
\medskip

\centerline{\bf $n=3$ case}
\medskip

For $n=3$ we have
$$
P=\alpha,\quad Q=\beta.
$$
Thus
$$
\aligned
f(x)&=(\alpha^2+2\alpha\cdot \beta)x + \beta^2x^3\\
\phi(x)&=\alpha^2+(2\alpha\cdot \beta + \beta^2)x^2.
\endaligned
$$
Now, one can easily get the initial data of the isofocal
deformations, which correspond to completely decomposable
situation for $n=3$.
\medskip
\begin{proposition}
The initial data of completely decomposable situation are given up
projective-linear transformations, by the formulae
$$
\aligned \alpha_1&=0\\
\alpha_{23}&=\pm \sqrt{-\frac{\alpha^2+2\alpha \beta}{\beta^2}}\\
m_1&=-\frac{\beta^2(2\alpha\beta +
\beta^2)}{\alpha^2+2\alpha\beta}\\
m_2&=m_3=\frac{\alpha^2(\alpha^2+2\alpha\beta)+
\beta^2(2\alpha\beta + \beta^2)}{2(\alpha^2+2\alpha\beta)}.
\endaligned
$$
\end{proposition}
\medskip

\centerline{\bf Case $n=5$}
\medskip

For $n=5$ we have
$$
P(x)=\alpha + \gamma x^2,\quad Q(x)=\beta.
$$
Thus,
$$
\aligned
f(x)&=x((\alpha^2+2\alpha\cdot \beta) + (\beta^2+2\alpha \gamma + 2\gamma\beta)x^2
+\gamma x^4)\\
\phi(x)&=\alpha^2+(2\alpha \beta +2\alpha \gamma+ \beta^2)x^2 +
(\gamma^2+2\beta\gamma)x^4.
\endaligned
$$
Denote by
$$
\aligned A&=\frac{-(2\alpha\gamma+2\gamma\beta
+\beta^2)+\sqrt{D}}{2\gamma^2}\\
B&=\frac{-(2\alpha\gamma+2\gamma\beta
+\beta^2)-\sqrt{D}}{2\gamma^2}\\
D&=2\alpha \gamma +
2\gamma\beta+\beta^2-4\gamma^2(\alpha^2+2\alpha\beta).
\endaligned
$$
Now we have
\medskip
\begin{proposition}
The initial data of completely decomposable situation are given up
projective-linear transformations, by the formulae
$$
\aligned \alpha_1&=0\\
\alpha_{23}&=\pm \sqrt {A}\\
\alpha_{45}&=\pm \sqrt {B}\\
m_1&=\frac{\gamma^2+2\alpha\gamma}{AB}\\
m_2&=m_3=\frac{1}{2(B-A)}(2\alpha\gamma+2\beta\alpha +
\beta^2+\frac{A+B}{AB}(\gamma^2+2\alpha\gamma) -
A\alpha^2+\frac{\gamma^2+2\alpha\gamma}{B})\\
m_4&=m_5=\frac{1}{2(B-A)}(2\alpha\gamma+2\beta\alpha +
\beta^2+\frac{A+B}{AB}(\gamma^2+2\alpha\gamma) -
B\alpha^2+\frac{\gamma^2+2\alpha\gamma}{A}).
\endaligned
$$
\end{proposition}
\medskip

\centerline{\bf Case $n=7$}
\medskip

For $n=7$ we have
$$
P=\alpha + \gamma x^2,\quad Q= \beta + \delta x^2,
$$
and
$$
\aligned \phi&=x(\alpha^2 + 2\alpha\beta + (2(\alpha \gamma +
\gamma \beta + \alpha \delta) x^2 + (\gamma^2 + 2(\gamma \delta +
\beta \delta)x^4 + \delta^2 x^6)\\
f&=\alpha^2 + (2(\alpha \gamma + \alpha \beta) + \beta ^2)x^2 +
(\gamma^2+2(\gamma \beta + \alpha \delta + \beta \delta))x^4 +
(2\gamma \delta + \delta^2)x^6.
\endaligned
$$
Denote by $A_1, A_2, A_3$ the three roots of the polynomial
$$
(\alpha^2 + 2\alpha\beta + (2(\alpha \gamma + \gamma \beta +
\alpha \delta) x + (\gamma^2 + 2(\gamma \delta + \beta \delta)x^2
+ \delta^2 x^3,
$$
and the seven zeros of the polynomial $\phi$ are
$$
\aligned
\alpha_0&=0,\quad \alpha_{12}=\pm \sqrt {A_1}\\
\alpha_{34}&=\pm \sqrt {A_2},\quad \alpha_{56}=\pm \sqrt {A_3}.
\endaligned
$$
To calculate the initial weights, introduce the new coordinates
$$
\aligned
Z_1 &= m_2 + m_3,\quad \hat Z_1 = m_2 - m_3\\
Z_2 &= m_4 + m_5,\quad \hat Z_2 = m_4 - m_5\\
Z_3 &= m_6 + m_7,\quad \hat Z_3 = m_6 - m_7
\endaligned
$$
Now, one gets the initial weights from the system of linear
equations:
$$
m_1=\frac{2\gamma \delta + \delta^2}{\alpha^2};
$$
$$
\aligned Z_1 + Z_2 + Z_3 &= d_1\\
\alpha_2^2Z_1 + \alpha_4^2Z_2 + \alpha_6^2Z_3 &= d_2\\
\alpha_2^2(\alpha_4^2+\alpha_6^2)Z_1 + \alpha_4^2(\alpha_2^2+\alpha_6^2)Z_2 + \alpha_6^2(\alpha_4^2+\alpha_2^2)Z_3 &= d_3\\
\alpha_2\hat Z_1 + \alpha_4\hat Z_2 + \alpha_6 \hat Z_3 &= 0\\
\alpha_2(\alpha_4^2+\alpha_6^2)\hat Z_1 + \alpha_4(\alpha_2^2+\alpha_6^2)\hat Z_2 + \alpha_6(\alpha_4^2+\alpha_2^2)\hat Z_3 &= 0\\
\alpha_2\alpha_4^2\alpha_6^2\hat Z_1 +
\alpha_4\alpha_2^2\alpha_6^2\hat Z_2 +
\alpha_6\alpha_4^2\alpha_2^2\hat Z_3 &= 0,
\endaligned
$$
where
$$
\aligned d_1&=\alpha^2-\frac{2\gamma\delta + \delta^2}{\alpha^2}\\
d_2&=(\gamma^2+2(\gamma \beta + \alpha \delta) + 2\beta
\delta)-\alpha^2(\gamma^2+2\gamma \delta + 2\beta \delta)\\
d_3&=(\gamma^2+2(\gamma \beta + \alpha \delta)-\frac{2\gamma
\delta + \delta^2}{\alpha^2}(2(\alpha\gamma + \gamma \beta +
\alpha \delta)+\beta^2).
\endaligned
$$
\medskip
\section{Concluding remarks}\label{sec:conclusion}
\medskip
\subsection{Bifocal transformations}
\medskip

After studying isofocal deformations, we may pass to a new class
of deformations, to {\it bifocal} deformations, defined by the
relation
\begin{equation}\label{eq:defbifocal}
B_t(z):=\frac{m_1(t)}{z-\alpha_1(t)}+\dots
+\frac{m_n(t)}{z-\alpha_n(t)}=\frac{f(z)+tg(z)}{\phi(z)+tf(z)},
\end{equation}
where polynomials $\phi$ and $f$ are the same as before, and $g$
is a new polynomial of degree $n-1$.

By applying the Moser trick once again and using the Flashka
coordinates, we come to the following
\medskip
\begin{proposition}\label{prop:bifocal}
The dynamical system (\ref{eq:defbifocal}) reduces in the Flashka
coordinates to  the form
\begin{equation}\label{eq:flashka}
\aligned \dot a_1=0 \quad \dot b_1&=0 \\
\dot a_2=0 \quad \dot b_2&=0 \\
\dots \\
\dot a_{n-2}=0 \quad b_{n-2}&=0\\
\dot a_{n-1}=\frac{b_{n-1}}{a_{n-1}} \quad \dot b_{n-1}&=-1\\
 \dot b_n&=-1.
 \endaligned
\end{equation}
\end{proposition}
\medskip

Denote zeros of the polynomial $G_t(z)=f(z)+tg(z)$ as
$f_1(t),\dots, f_{n-1}(t)$. Here $f_1(0)=f_1,\dots,
f_{n-1}(0)=f_{n-1}$ are zeros of the polynomial $f$. In bifocal
case, focal points of Marden curves evolve during the time as
zeros of $G_t$.

\medskip

Denote by $a(\alpha):=a(\alpha_1,\dots, \alpha_n)$ the matrix
\begin{equation}
a(\alpha)=\left(\begin{array}{cccc} 1 & 1 &
\dots & 1\\
\sigma^1_1(\alpha_1,\dots, \alpha_n) & \sigma^2_1(\alpha_1,\dots,
\alpha_n) & \dots & \sigma^n_1(\alpha_1,\dots, \alpha_n)\\
\dots & \dots & \dots & \dots\\
\sigma^1_{n-1}(\alpha_1,\dots, \alpha_n) &
\sigma^2_{n-1}(\alpha_1,\dots, \alpha_n) & \dots &
\sigma^n_{n-1}(\alpha_1,\dots, \alpha_n)
\end{array}\right)
\end{equation}
\medskip
Denote also the columns
\medskip
$$
\aligned \dot \alpha &= \left(\begin{array}{cccc} \dot \alpha_1 &
\dot \alpha_2 & \dots & \dot \alpha_n
\end{array}\right)^T\\
\overrightarrow B &= \left(\begin{array}{cccc}B_1 & B_2 & \dots &
B_{n}
\end{array}\right)^T\\
\overrightarrow m &= \left(\begin{array}{cccc}m_1 & m_2 & \dots &
m_{n}
\end{array}\right)^T\\
\dot m&=\left(\begin{array}{cccc}\dot m_1 & \dot m_2 & \dots &
\dot m_{n}
\end{array}\right)^T
\endaligned
$$
\medskip
\begin{proposition}\label{prop:eqmotion}
The first part of the differential equations for isofocal and
bifocal deformations may be written down in the following way:
\begin{equation}\label{eq:finaleqmotion1}
\dot \alpha= a^{-1}(\alpha_1,\dots, \alpha_n) \overrightarrow B.
\end{equation}
\medskip
The second part of the equations of motion are equivalent to
\begin{equation}\label{eq:finaleqmotion2}
\dot m = a^{-1}(\alpha)(\overrightarrow C-\dot a(\alpha)
\overrightarrow m),
\end{equation}
where $\dot a(\alpha):= \dot a(\alpha_1,\dots, \alpha_n)$ denotes
the matrix obtained from $a(\alpha_1,\dots, \alpha_n)$ by
differentiating each matrix entry of the matrix $a$ with respect
to the time, taking into the account the equations
(\ref{eq:finaleqmotion1}). The column $\overrightarrow C$ is equal
to the differential of the column $\overrightarrow B$ with respect
to the time. In the case of isofocal deformations,
$\overrightarrow C$ is equal to the zero-vector.
\end{proposition}
\medskip
The second part of the equations of motion follows by
differentiation with respect to time, from the relation,
$$
a(\alpha_1,\dots, \alpha_n)\overrightarrow m = \overrightarrow B.
$$
\medskip
\subsection{Toma-Trautmann case of a conic component}
\medskip
There is one more interesting case of conic component of a
Darboux-Poncelet curve. Such cases have been studied by Trautmann
systematically in \cite {Tr1}, following an example constructed by
Toma. It is the case where a conic component $C$ of a
Poncelet-Darboux curve $S$ of degree $n-1$ is by itself a
Poncelet-Darboux curve. In other words, in this case the conic $C$
is $3$-Ponclet related to the caustic $K$. Suppose that $S$ is
defined by a pencil of polynomials of degree $n$ generated by $f$
and $\phi$. Applying the same arguments as before we conclude that
there should exist four parameters $t_1, t_2, t_3, t_4$ such that
\begin{equation}\label{eq:decomp3}
\aligned \Phi_{t_1}(z):=\phi(z)+t_1f(z)&=(z-a_1)Q_1^2(z)S_1(z)\\
\Phi_{t_2}(z):=\phi(z)+t_2f(z)&=(z-a_2)Q_2^2(z)S_2(z)\\
\Phi_{t_3}(z):=\phi(z)+t_3f(z)&=(z-a_3)Q_3^2(z)S_3(z)\\
\Phi_{t_4}(z):=\phi(z)+t_4f(z)&=(z-a_4)Q_4^2(z)S_4(z);
\endaligned
\end{equation}
where polynomials $Q_i$ are of degree one and polynomials $S_i$
are of degree $n-3$. Thus we get
\medskip
\begin{proposition}\label{prop:tomatraut}
Denote by $\Gamma$ hyperelliptic curve of genus $g=2n-5$ defined
by the equation
$$
y^2=(z-a_1)(z-a_2)(z-a_3)(z-a_4)S_1(z)S_2(z)S_3(z)S_4(z),
$$
and by $\Gamma_2$ the elliptic curve
$$
Y^2=(X+t_1)(X+t_2)(X+t_3)(X+t_4).
$$
The Toma-Trautmann case of conic component of a Poncelet-Darboux
curve of degree $n-1$ is related to the above covering of the
hyperelliptic $\Gamma$ of genus $2n-5$ over the elliptic curve
$\Gamma_2$.
\end{proposition}
\medskip
Although the hyperelliptic coverings of elliptic curves are much
more complicated than the elliptic coverings of elliptic curves
discussed above, an intensive   study of former coverings has been
done, see for example \cite {AP}, \cite {Tre} and references
therein. It sounds as an interesting question to describe
hyperelliptic coverings which arise in the Toma-Trautmann case.
\medskip
\subsection{Infrapolynomials and the positivity problem in isofocal
dynamics}
\medskip
From \cite{Ma} we know that one important interpretation of the
Marden Theorem is connected with the study of {\it
infrapolynomials} (see Section 5 of Chapter 1 of \cite {Ma} for
definitions). After Fekete it is known that for given closed
bounded set $E$ containing at least $n+1$ points, a polynomial $p$
of degree $n$ is an infrapolynomial if there exist an integer $k$,
such that $n \le k \le 2n$, a set of {\bf positive} constants
$m_j$, such that $m_0+\dots+m_k=1$ and a set of $k+1$ points
$\{z_0, z_1,\dots, z_k\}\subset E$ such that $p(z)$ is a factor of
the polynomial $F(z)$:
$$
F(z)=\Omega(z)\sum_{j=0}^k\frac{\lambda_j}{z-z_j},\quad
\Omega(z)=\prod_{j=0}^k(z-z_j).
$$
The last theorem of Fekete and the connection with the study of
infraplynomials thus motivate a study of positivity conditions in
our isofocal dynamics: 1) to describe intervals in isofocal
dynamics where all weights are positive; 2) to describe the
initial values for which the weights are positive during entire
evolution; 3) to relate more closely isofocal dynamics and the
theory of infrapolynomials.

\medskip
\section{Appendix}
\medskip

{\bf The case $n=3$.} Even in the simplest case $n=3$ when the
curve $C$ is a conic, inscribed in the triangle formed by the
zeros of a polynomial of degree 3 and tangent to the sides of the
triangle at their midpoints, the result of the Siebeck theorem is
nontrivial and interesting and attracted lot of attention not only
in the past but also nowadays (see for example \cite{Ka}). We are
going to present main points of the proof of the Siebeck theorem
for $n=3$ following mostly \cite{Ka} and references therein.

We start with well known focal properties of ellipses and caustic
properties of elliptical billiards. They are going to play central
role in the proof of the Siebeck theorem.
\medskip
\begin{lemma} [Focal property of the ellipse]%
\label{prop:focal.property}%
Let $\mathcal E$ be an ellipse with foci $F_1,F_2$ and
$A\in\mathcal E$ an arbitrary point. Then segments $AF_1$, $AF_2$
satisfy the billiard law on $\mathcal E$.
\end{lemma}
\medskip
\begin{lemma}\label{prop:ellipse.caustic}
Let two lines satisfy the billiard law on the ellipse $\mathcal
E$. If one of the lines is tangent to the ellipse $\mathcal E'$
that is confocal with $\mathcal E$, then the other one is also
tangent to $\mathcal E'$.
\end{lemma}
\medskip
\begin{corollary}
Let $\mathcal T$ be a billiard trajectory within ellipse $\mathcal
T$. If $\mathcal C$ is a conic confocal to $\mathcal E$ such that
one segment of $\mathcal T$ is tangent to $\mathcal C$, then all
segments of $\mathcal T$ are tangent to $\mathcal C$.
\end{corollary}
\medskip
\begin{corollary}
Let $B$ be a point outside the ellipse $\mathcal E'$ with focal
points $F_1$ and $F_2$. Denote tangents to the ellipse from the
point $B$ as $BB_1$ and $BB_2$, where $B_i$ are points of contact
with the ellipse. Then the angles $B_1BF_1$ and $B_2Bf_2$ are
equal.
\end{corollary}
\medskip
The first observation in the proof is that the statement is
invariant to applications of affine transformations. Now, we
assume that a real plane is identified with the field of complex
numbers and that numbers $\alpha_1, \alpha_2, \alpha_3$ correspond
to numbers $-1, 1, w$ where $w$ is in upper half-plane.

The polynomial $P$ gets the form
$$
P(z)=(z-1)(z+1)(z-w)=z^3-wz^2-z+w,
$$
and the derivative is
$$
DP(z)=3\left(z^2-\frac{2}{3}wz-\frac{1}{3}\right).
$$
Denote the zeros of the derivative as $z_4, z_5$. After some
analysis of the formulae
$$
\aligned z_4+z_5&=\frac{2}{3}w\\
z_4z_5&=-\frac{1}{3}
\endaligned
$$
we conclude that both of the points $z_4$ and $z_5$ are in the
upper half plane. Then, from the second formula we conclude that
the sum of their arguments $\theta_4$ and $\theta_5$ is equal to
$\pi$. Denote the line  connecting $z_4$ and the origin $O$ as
$L_{z_4O}$ and the line  connecting $z_5$ and the origin $O$ as
$L_{z_5O}$. As a result, we come to the conclusion that the angle
between the line $L_{z_4O}$ and negative real semi-axis is equal
to the angle of the line $L_{z_5O}$ and the positive real
semi-axis.  By application of focal properties of ellipses,
Proposition (\ref{prop:focal.property}) and Proposition
(\ref{prop:ellipse.caustic}), we see that there is an ellipse
$E_1$, with $z_4, z_5$ as foci, such that it touches the real axis
at the origin. In other words, it touches the segment $[\alpha_1,
\alpha_2]$ at the midpoint.

By symmetry, we see that there are confocal ellipses $E_1$, $E_2$
and $E_3$ with focal points at the zeros of derivative of the
polynomial, each of which touches one of the sides of the triangle
in the midpoint.

Next, we need to show that these three ellipses coincide. Again,
we apply affine transformation and transform the triangle to a new
one with vertices $0, 1, w$, where $w$ is again in the upper
half-plane.

The polynomial $P$ gets the form
$$
P(z)=z(z-1)(z-w)=z^3-(1+w)z^2+wz,
$$
and the zeros $z_4, z_5$ of the derivative
$$
DP(z)=3\left(z^2-\frac{2}{3}(1+w)z+w\right)
$$
satisfy
$$
\aligned z_4+z_5&=\frac{2}{3}(1+w)\\
z_4z_5&=\frac{w}{3}.
\endaligned
$$
After analysis of the last formulae we again conclude that $z_4,
z_5$ are in the upper half-plane and that the sum of their
arguments is equal to the argument of $w$. As a result, this time
we see that the angle between the line $L_{z_4O}$ and the line
$L_{Ow}$ is equal to the angle of the line $L_{z_5O}$ and the
positive real semi-axis. From previous considerations, we know
that the ellipse $E_1$ with $z_4, z_5$ as focal points touches
real axis at the point $1/2$. One can easily see that the origin
$O$ is outside the ellipse $E_1$. Now, we apply another well known
focal property of ellipses (see Corollary of Proposition
(\ref{prop:ellipse.caustic}) to the ellipse $E_1$, point $O$
outside $E_1$ and two tangents $t_1$ and $t_2$ to $E_1$ from the
point $O$. The angle between $t_1$ and $L_{Oz_5}$ is equal to the
angle between $t_2$ and $L_{Oz_4}$. One tangent,  $t_1$ is the
real axis. Thus, the second tangent coincides with the line
$L_{Ow}$. This proves that the line $L_{Ow}$ is tangent to the
ellipse $E_1$. But, again by focal properties of ellipses, (see
Proposition (\ref{prop:focal.property})) among confocal ellipses,
there is only one tangent to a given line. This shows that
$E_1=E_2$ and $E_1$ touches the segment $[O,w]$ at the midpoint.
This finishes the proof of the Siebeck theorem in the case $n=3$.

The Lemmae and Corollaries from the beginning of the Appendix play
very important role in the study of billiard systems within
ellipses. These systems provide mechanical interpretation of the
Poncelet theorem and this is the way how potential connection of
Marden's and Poncelet's theorem has been anticipated.

\medskip
\subsection*{Acknowledgements}

The research was partially supported by the Serbian Ministry of
Science and Technology, Project {\it Geometry and Topology of
Manifolds and Integrable Dynamical Systems}.  The last part of the
paper has been written during a visit to the IHES. The author uses
the opportunity to thank the IHES for hospitality and outstanding
working conditions.

\newpage\thispagestyle{empty}


\vspace*{20mm} 

\begin{thebibliography}{99}

\bibitem{AP}
R. D. M. Accola, E. Previato, {\it Covers of tori: genus two}
Lett. Math. Phys. {\bf 76} (2006) 135-161


\bibitem{BM}
W. Barth, J. Michel, {\it Modular curves and Poncelet polygons}
Math. Ann. {\bf 295} (1993) 25-49


\bibitem{Ber}
M. Berger, {\it Geometry}, Springer-Verlag, Berlin, 1987.


\bibitem{Dar1}
G. Darboux, {\it Principes de g\'eom\'etrie analytique},
Gauthier-Villars, Paris (1917) 519 p.

\bibitem{Dar3}
G. Darboux, {\it  Le\c{c}ons sur la th\'eorie g\'en\'erale des
surfaces et les applications g\'eo\-m\'etri\-ques du calcul
infinitesimal}, volumes 2 and 3, Gauthier-Villars, Paris, 1887,
1889.


\bibitem{DR}
V. Dragovi\'c, M. Radnovi\'c, {\it Geometry of integrable
billiards and pencils of quadrics}, Journal Math. Pures Appl. {\bf
85} (2006), 758-790.
\newline
arXiv: math-ph/0512049

\bibitem{DR2}
V. Dragovi\'c, M. Radnovi\'c, {\it Hyperelliptic Jacobians as
Billiard Algebra of Pencils of Quadrics: Beyond Poncelet Porisms},
Adv.  Math., {\bf 219} (2008).
\newline
arXiv: math-ph/0710.3656

\bibitem{Fl}
H. Flashka, {\it The Toda Lattice} I, Phys. Rev. B {\bf 9}, (1974)
1924-1925

\bibitem{Ha}
R. Hartshorne,{\it Stable vector bundles of rank 2 on $\mathbb
{P}_3$}, Math. Ann. {\bf 238} (1978), 229-280.

\bibitem{HN}
A. Hirschowitz, M. S. Narasimhan, {\it Fibres de 'tHooft sp\'
eciaux et applications}, Proc. Nice Conf. 1981, 142-163,
Birkhauser (1982)


\bibitem{Jac2}
C. Jacobi, {\it Fundamenta nova theoriae functiorum ellipticarum}
(1829)

\bibitem{Ka}
D. Kalman, {\it An Elementary Proof of Marsden's Theorem} The
American Mathematical Monthly, {\bf 115} (2008) 330-337

\bibitem{Ma}
M. Marden, {\it Geometry of Polynomials}, AMS, Math. Surveys, No
3, Second edition (1966) 243 p.


\bibitem{Mo}
J. Moser, {\it Finitely many mass points on the line under the
influence of an exponential potential--an integrable system},
Lecture Notes in Physics, {\bf 28} Springer (1975), 467-497

\bibitem{Mo1}
J. Moser, {\it Three integrable Hamiltonian systems connected with
isospectral deformations}, Adv. Math., {\bf 16} (1975) 197-220


\bibitem{NT}
M. S. Narasimhan, G. Trautmann, {\it Compactification of
$M_{P_3}(0,2)$ and Poncelet pairs of conics}, Pacific J. of Math.
{\bf 145} (1990) 255-365

\bibitem{Pon}
J. V. Poncelet, {\it Trait\'e des propri\'et\'es projectives des
figures}, Mett-Paris, 1822.

\bibitem{Si}
J. Siebeck, {\it Ueber eine neue analytische Behandlungweise der
Brennpunkte}, J. Reine Angew. Math. {\bf 64} (1864) 175-182

\bibitem{Tr}
G. Trautmann, {\it Poncelet curves and theta characteristics}
Expositiones Mathematicae, {\bf 6} (1988) 29-64

\bibitem{Tr1}
G. Trautmann, {\it Decomposition of Poncelet curves and instanton
bundles}, A. St. Univ. Ovidius Constanta, {\bf 5} (1997) 105-110

\bibitem{Tre}
A. Treibich, {\it Hyperelliptic tangential covers and finite-gap
potentials} Russian Math. Survays {\bf 56} (2001) 1107- 1151


\end{thebibliography}
\end{document}